\documentclass[12pt,reqno]{amsart}

\usepackage[T1]{fontenc}
\usepackage[utf8]{inputenc}
\usepackage{lmodern}
\usepackage{amsmath,amssymb,amsthm,mathrsfs}
\usepackage{microtype}
\usepackage{enumitem}
\usepackage[hidelinks]{hyperref}

\allowdisplaybreaks[2]

\theoremstyle{plain}
\newtheorem{theo}{Theorem}[section]
\newtheorem{coro}[theo]{Corollary}
\newtheorem{lemm}[theo]{Lemma}
\newtheorem{prop}[theo]{Proposition}
\theoremstyle{definition}
\newtheorem{defi}[theo]{Definition}
\theoremstyle{remark}
\newtheorem{rem}[theo]{Remark}

\numberwithin{equation}{section}

\title[Li--Yau estimates for a Finsler logarithmic equation] {A Li--Yau gradient estimate for the Finslerian logarithmic Schr\"odinger equation} \author{Duoxin Li, Zisu Zhao$\dagger$}
\address{School of Mathematical Sciences, Huaqiao University, Quanzhou, Fujian, 362021, China.}
\email{duoxinli52math@163.com}
\address{Institute of Mathematics, Academy of Mathematics and Systems Science, Chinese Academy of Sciences, Zhongguancun East Road 55, Beijing 100190,
China} \email{steven1998zhao@outlook.com} 
\thanks{$\dagger$The first author is supported by the Scientific Research Funds of Huaqiao University (Grant No. 605-50Y25056). The second author was partially supported by the China Scholarship Council (Grant No.\ 202306310146). The authors would like to thank Professor Chunhui Qiu for his consistent encouragement and helpful suggestions.} 

\begin{document}

\begin{abstract}
We establish a local Li--Yau type gradient estimate for a nonlinear parabolic equation with logarithmic reaction on a forward complete Finsler metric measure manifold. The
coefficients may depend on both space and time. The estimate uses a weighted Ricci curvature lower bound and explicit cutoff bounds for the divergence-form linearized Laplacian. We explain the role of
this localization assumption when the reference direction varies with the solution. Applications include a local Harnack inequality for the equation with a
potential and an upper bound for positive stationary solutions under a suitable cutoff exhaustion. An appendix gives a variational existence proof for positive bounded initial
data, with quantitative bounds on every finite time interval, which is new as far as the authors know.
\end{abstract}

\maketitle

\section{Introduction}

We study positive solutions of the parabolic equation
\begin{equation}\label{1.1}
    u_t=\Delta u+a(x,t)u\log u+b(x,t)u
\end{equation}
on a Finsler metric measure manifold $(M,F,\mu)$. Here $\Delta$ denotes the nonlinear Finsler Laplacian associated with $F$ and the smooth positive measure $\mu$. Equation~\eqref{1.1} is the real
parabolic equation associated with the logarithmic Schr\"odinger equation. Logarithmic nonlinearities occur in models related to wave mechanics, information theory, and transport
\cite{BJD1991,MSD2003,YK1978}. Our concern is the effect of curvature and the Finsler structure on the behavior of positive solutions.

Gradient estimates for logarithmic equations have been studied in several geometric settings. Ma obtained estimates for a corresponding elliptic equation
on complete noncompact Riemannian manifolds \cite{LM2006}. Huang and Ma treated the parabolic equation with constant coefficients under a weighted Ricci
curvature lower bound \cite{HM2010}, while Cao and Zhang studied differential Harnack estimates for parabolic equations in geometric settings \cite{CXZZ2011}. Dung and Khanh established
Hamilton--Souplet--Zhang type estimates for a general heat equation \cite{DK2015}, following the method of Souplet and Zhang \cite{SPZQS2006}. For coefficients depending on space
and time, Wang proved estimates for elliptic and parabolic logarithmic equations \cite{WangJ2023}. Related existence questions on locally finite graphs were studied by Chang, Wang, and Yan, who obtained ground
states for a stationary logarithmic Schr\"odinger equation \cite{CWY2023}.

The Finsler setting introduces a further nonlinearity through the gradient itself. Ohta and Sturm constructed the Finsler heat flow as a gradient flow of a convex energy on
$L^2(M,\mu)$ \cite{OhtaSturm2009}. They later established the Bochner--Weitzenb\"ock formula and Li--Yau estimates \cite{OSSKT2014}. Xia obtained further Li--Yau estimates on compact and complete
noncompact Finsler manifolds \cite{XiaQ2022}. Shen introduced mixed weighted Ricci curvature and investigated Laplacian comparison and gradient estimates for Finslerian
equations \cite{ShenB2023}. Shen and Zhu also studied a nonlinear parabolic equation on Finsler metric measure spaces \cite{ShenB2025}.

The argument below follows Wang's choice of auxiliary function and uses the integrated Finsler Bochner inequality. Two points require care. The operators $\nabla$ and $\Delta$ are nonlinear, whereas the operators
$\nabla^{\nabla u}$ and $\Delta^{\nabla u}$ are linear in their argument once the reference field is fixed. Moreover, $g_{\nabla u}$ need not extend smoothly across the critical set of $u$.
We therefore keep the linearized operators in divergence form and use their weak formulations throughout.

Localization also requires control of the actual reference field. The divergence of $g^{ij}(x,\nabla u)\phi_j$ contains derivatives of $\nabla u$ in addition to
derivatives of the metric in the base variable. We state the required cutoff bounds explicitly in the main theorem. Section~\ref{subsec:cutoffs} explains how a valid radial comparison
yields these bounds and why derivatives of the reference field must be retained. In particular, the cutoff hypothesis is an additional analytic assumption and is not asserted to follow from the
mixed weighted Ricci bound alone.

For a function $h$ on $B_p(2R)\times(0,T]$, write
\[
    [h]^+=\operatorname*{ess\,sup}_{B_p(2R)\times(0,T]} \max\{h,0\}.
\]
For a constant $s$, write $s^+=\max\{s,0\}$. We use $F_V(X)=\sqrt{g_V(X,X)}$ for the norm induced by the fundamental tensor with reference vector $V$. The geometric
and analytic definitions used in the following statement are recalled in Section~\ref{sec:preliminaries}.

\begin{theo}\label{MainTheorem}
Let $(M,F,\mu)$ be a forward complete $n$-dimensional Finsler metric measure manifold, and let $N>n$ be finite. Fix $p\in M$, $R>0$, and $T>0$, with $\overline{B_p(2R)}$
disjoint from $\partial M$. Assume that the misalignment on $B_p(2R)$ is at most $\alpha<\infty$, and that the mixed weighted Ricci curvature satisfies
\[
    {}^m\operatorname{Ric}^N_W(V)\ge-KF^2(V)
\]
for all nonzero $V,W$ based in $B_p(2R)$, where $K\ge0$.

Let $u$ be a positive solution of \eqref{1.1} on a neighborhood of $\overline{B_p(2R)}\times(0,T]$, and suppose $u\le D$ for a constant $D>0$. Assume the interior regularity in \eqref{eq:solution-regularity}. Let $a$ and $b$ be
$C^2$ in space and $C^1$ in time, and suppose their weak linearized Laplacians have locally integrable representatives. Suppose that $a$, $b$, $a_t$, $F_{\nabla u}(\nabla^{\nabla u}a)$, and $F_{\nabla u}(\nabla^{\nabla
u}b)$ are bounded, $\Delta^{\nabla u}b$ is bounded below, and $[\Delta^{\nabla u}a+a_t]^+<\infty$. All bounds are taken on $B_p(2R)\times(0,T]$. The linearized Laplacians and their bounds are understood weakly, including across
$\{du=0\}$.

Assume, in addition, that there is a locally Lipschitz spatial cutoff $\phi$, equal to one on $B_p(R)$ and compactly supported in $B_p(2R)$, with $0\le\phi\le1$, such that
\begin{equation}\label{eq:cutoff-hypotheses}
    \frac{F_{\nabla u}^2(\nabla^{\nabla u}\phi)}{\phi} \le\frac{\alpha C_1^2}{R^2},\qquad \frac{2F_{\nabla u}^2(\nabla^{\nabla u}\phi)}{\phi} -\Delta^{\nabla u}\phi\le B
\end{equation}
on $\{\phi>0\}\times(0,T]$, for constants $C_1>0$ and $B\ge0$. These inequalities are also understood weakly.

Choose $A>[a]^+$ and $E\ge0$ satisfying \eqref{eq:E-positive}--\eqref{eq:E-condition}. Such an $E$ depends only on $A$, $D$, $N$, and the coefficient bounds listed
above. Then
\begin{align}\label{1.2}
    &\frac{F^2(\nabla u)}{u^2} +(A+a)\log\frac{u}{D}-2\frac{u_t}{u} \notag\\
    &\quad\le4N\left\{ \frac1t+[a]^++B+\frac{N\alpha C_1^2}{R^2} \right.\notag\\
    &\qquad\left. +\frac{(2K+2+|\log D|-A)^+}{4} +\frac{2[E-a\log D]^+}{N} +\frac{[\Delta^{\nabla u}a+a_t]^+}{4(A-[a]^+)} \right\}
\end{align}
on $B_p(R)\times(0,T]$.
\end{theo}

The proof only uses the ordinary weighted Ricci lower bound obtained by taking $V=W$ in the mixed curvature assumption. Thus the mixed bound can be replaced by $\operatorname{Ric}^N\geq-K$, provided the cutoff hypothesis is
retained. When $F$ is Riemannian, the linearized operator is independent of $u$, and the usual weighted Laplacian comparison supplies these cutoffs.

Section~\ref{sec:preliminaries} recalls the analytic and geometric identities and discusses the localization assumption. Section~\ref{sec:gradient-proof} proves the local estimate and a separate
estimate for the equation with $a=0$. Section~\ref{sec:applications} gives the Harnack inequality and the stationary bound. Appendix~\ref{app:existence} constructs positive solutions and establishes their interior
regularity.

\section{Preliminaries}\label{sec:preliminaries}

We recall the definitions and identities needed in the proof. General references for Finsler geometry include \cite{AP1,DCS2012,ShenZ2001}; the analytic conventions used here follow \cite{Ohta2017,Ohta2021,OhtaSturm2009,OSSKT2014}.

\subsection{The Finsler gradient and the linearized Laplacian}

A smooth Finsler metric is a continuous function $F:TM\to[0,\infty)$ that is smooth and positive away from the zero section, satisfies $F(x,\lambda y)=\lambda F(x,y)$ for $\lambda>0$, and has
positive definite fundamental tensor
\[
    g_{ij}(x,y)=\frac12 \frac{\partial^2 F^2}{\partial y^i\partial y^j}(x,y), \qquad y\ne0.
\]
Reversibility is not assumed, so $F(x,-y)$ may differ from $F(x,y)$. For $V\ne0$, we write
\[
    \langle X,Y\rangle_V=g_V(X,Y),\qquad F_V(X)=\sqrt{g_V(X,X)}.
\]
In particular, $g_y(y,y)=F^2(y)$. The dual norm and the Legendre map are given by
\[
    F^*(x,\xi)=\sup_{F(x,y)=1}\xi(y),\qquad \mathcal L(y)=g_y(y,\cdot),\qquad \mathcal L(0)=0.
\]
The Legendre map is generally nonlinear. It is a homeomorphism on each fiber and a smooth diffeomorphism away from the zero section. The gradient of $f$ is $\nabla f=\mathcal L^{-1}(df)$, with $\nabla f=0$ where $df=0$. Thus
\[
    df(\nabla f)=F^2(\nabla f)=F^{*2}(df).
\]
For a nonzero reference field $V$, the linearized gradient is
\[
    \nabla^V f=g^{ij}(x,V)f_j\,\partial_{x^i}.
\]
It agrees with $\nabla f$ when $V=\nabla f$, but it need not agree with the nonlinear gradient for other choices of $V$.

Write $d\mu=\sigma(x)\,dx$ in a coordinate chart. The divergence, Laplacian, and linearized Laplacian are
\begin{equation}\label{eq:laplacians}
    \operatorname{div}_\mu X=\frac1\sigma\partial_i(\sigma X^i), \qquad \Delta f=\operatorname{div}_\mu(\nabla f), \qquad \Delta^V f=\frac1\sigma\partial_i(\sigma g^{ij}(x,V)f_j).
\end{equation}
Repeated coordinate indices are summed. In weak form,
\[
    \int_M\varphi\Delta f\,d\mu =-\int_M d\varphi(\nabla f)\,d\mu
\]
for $\varphi\in C_c^\infty(M)$. Similarly,
\[
    \int_M\varphi\Delta^V f\,d\mu =-\int_M g_V^*(d\varphi,df)\,d\mu.
\]
The operators with superscript $V$ are linear in $f$ when $V$ is fixed, and satisfy the usual product rule
\begin{equation}\label{eq:product-rule}
    \Delta^V(fh)=f\Delta^Vh+h\Delta^Vf +2g_V^*(df,dh).
\end{equation}
For $h\in C^2(\mathbb R)$ and a sufficiently regular $f$, the corresponding chain rule is
\[
    \Delta^V(h\circ f) =h'(f)\Delta^Vf+h''(f)g_V^*(df,df).
\]

We denote the noncritical set of $u$ by $M_u=\{x\in M\,\vert\,du(x)\ne0\}$. To use the linearized operator on all of $M$, choose a measurable nonzero
reference field agreeing with $\nabla u$ on $M_u$ and use its fundamental tensor on the critical set. The resulting coefficients are uniformly elliptic on compact coordinate
domains. Identities involving this extension are understood weakly, as in \cite{OSSKT2014}. Bounds on the linearized Laplacians of the coefficients and
of the cutoff refer to this same choice of extension.

The misalignment on a set $\Omega\subset M$ is
\begin{equation}\label{eq:misalignment}
    \alpha_\Omega =\sup_{x\in\Omega}\ \sup_{V,W,Z\in T_xM\setminus\{0\}} \frac{g_V(Z,Z)}{g_W(Z,Z)}.
\end{equation}
This quantity compares the fundamental tensors associated with different reference directions \cite{ShenB2023}. It is at least one and is finite on every compact
subset. A bound $\alpha_\Omega\leq\alpha$ gives
\[
    \alpha^{-1}g_W\leq g_V\leq\alpha g_W, \qquad \alpha^{-1}g_W^*\leq g_V^*\leq\alpha g_W^*.
\]
These comparisons are related to the uniform convexity and smoothness estimates discussed in \cite{BCL1994,Ohta2017}.

\subsection{Weighted curvature and the Bochner inequality}

Let $D^y$ denote the Chern connection with reference direction $y$. On $M_f$, the Hessian of $f$ is the symmetric bilinear form
\[
    \nabla^2f(X,Y) =g_{\nabla f}(D_X^{\nabla f}\nabla f,Y).
\]
Its coordinate expression is $f_{ij}-\Gamma^k_{ij}(x,\nabla f)f_k$. The distortion and the $S$-curvature are
\[
    \tau(x,y)=\log\frac{\sqrt{\det g_{ij}(x,y)}}{\sigma(x)}, \qquad S(y)=\left.\frac{d}{ds} \tau(\gamma(s),\dot\gamma(s))\right|_{s=0},
\]
where $\gamma$ is the geodesic with initial velocity $y$. We write $\dot S(y)=\left.\frac{d}{ds}S(\dot\gamma(s))\right|_{s=0}$. For $n<N<\infty$, the weighted Ricci curvature is
\begin{equation}\label{eq:weighted-Ricci}
    \operatorname{Ric}^N(y) =\operatorname{Ric}(y)+\dot S(y) -\frac{S(y)^2}{N-n}.
\end{equation}
All terms in this formula are homogeneous of degree two in $y$. Thus $\operatorname{Ric}^N\geq-K$ means $\operatorname{Ric}^N(y)\geq-KF^2(y)$ for every $y$. The familiar identity between the Hessian and the Laplacian is
\[
    \Delta f=\operatorname{tr}_{g_{\nabla f}}\nabla^2f-S(\nabla f).
\]

For completeness, let $R_y(X,y)y$ be the curvature endomorphism in the Jacobi equation. Shen's mixed weighted Ricci curvature can be written as
\begin{equation}\label{eq:mixed-Ricci}
    {}^m\operatorname{Ric}^N_W(y) =g^{ij}(x,W)g_y(R_y(e_i,y)y,e_j) +\dot S(y)-\frac{S(y)^2}{N-n},
\end{equation}
where $\{e_i\}$ is any coordinate basis and $W,y\ne0$ \cite{ShenB2023}. In particular, ${}^m\operatorname{Ric}^N_y(y)=\operatorname{Ric}^N(y)$. The mixed lower bound in the main theorem is
imposed for every pair of nonzero directions and therefore implies the ordinary weighted Ricci lower bound needed in the Bochner inequality.

If $f\in H^2_{\mathrm{loc}}\cap C^1$ and $\Delta f\in H^1_{\mathrm{loc}}$, the integrated Bochner inequality of Ohta and Sturm gives
\begin{equation}\label{eq:Bochner-prelim}
    \frac12\Delta^{\nabla f}F^2(\nabla f) -d(\Delta f)(\nabla f) \geq\operatorname{Ric}^N(\nabla f) +\frac{(\Delta f)^2}{N}
\end{equation}
in the weak sense \cite[Theorem~3.6]{OSSKT2014}. On the noncritical set, the corresponding identity contains $\operatorname{Ric}^{\infty}(\nabla f) +\|\nabla^2f\|_{\mathrm{HS}(g_{\nabla f})}^2$ on the
right. The weighted Cauchy--Schwarz inequality yields \eqref{eq:Bochner-prelim}. Its integrated form permits its use across the critical set.

\subsection{Cutoffs for the divergence linearization}
\label{subsec:cutoffs}

Let $r(x)=d_F(p,x)$ be the forward distance from $p$ and let $B_p(R)=\{x\in M\,\vert\,r(x)<R\}$. Forward completeness implies compactness of closed forward balls. Choose a smooth nonincreasing
function $\eta:[0,\infty)\to[0,1]$ with
\[
    \eta=1\ \text{on }[0,1],\qquad \eta=0\ \text{near }[2,\infty),\qquad \frac{(\eta')^2}{\eta}\leq C_1^2, \qquad \eta''\geq-C_2
\]
where $\eta>0$. Set $\phi=\eta(r/R)$. At regular points of $r$, $F^*(dr)=1$ and the misalignment bound implies
\begin{equation}\label{eq:radial-gradient-cutoff}
    \frac{F_V^2(\nabla^V\phi)}{\phi} \leq\frac{\alpha C_1^2}{R^2}.
\end{equation}
If the actual divergence operator with reference field $V$ satisfies
\begin{equation}\label{eq:radial-comparison-assumption}
    \Delta^V r \leq C(N,\alpha)\sqrt{\frac K{C(N,\alpha)}} \coth\left(r\sqrt{\frac K{C(N,\alpha)}}\right)+C_0, \qquad C(N,\alpha)=N+(\alpha-1)n-\alpha,
\end{equation}
in a sense that permits the distance cutoff, then
\begin{align}\label{eq:radial-B}
    B={}&\frac{2\alpha C_1^2+\alpha C_2}{R^2}\notag\\
    &+\frac{C_1}{R}\left\{ C(N,\alpha)\sqrt{\frac K{C(N,\alpha)}} \coth\left(R\sqrt{\frac K{C(N,\alpha)}}\right)+C_0 \right\}
\end{align}
is an admissible constant in \eqref{eq:cutoff-hypotheses}. For $K=0$, the expression involving $\coth$ is interpreted as $C(N,\alpha)/R$. Indeed,
\[
    \Delta^V\phi =\frac{\eta'}R\Delta^V r +\frac{\eta''}{R^2}g_V^*(dr,dr),
\]
and $\eta'$ is supported where $R\leq r\leq2R$. These formulas give \eqref{eq:radial-B} directly. At cut points one must use a barrier or a weak approximation for
which the same bounds hold; the usual distance support argument is due to Calabi \cite{Calabi1958}.

Comparison formulas of the form \eqref{eq:radial-comparison-assumption} motivate this choice of cutoff \cite{ShenB2023}. For a varying field $V$, however, the full divergence in
\eqref{eq:laplacians} must be retained. In coordinates,
\begin{equation}\label{eq:reference-derivative}
    \partial_i\bigl(g^{ij}(x,V(x))\bigr) =\left.\partial_{x^i}g^{ij}(x,y)\right|_{y=V(x)} +\left.\partial_{y^k}g^{ij}(x,y)\right|_{y=V(x)} \partial_i V^k.
\end{equation}
The second term does not vanish in general, even for a non-Riemannian Minkowski metric. A bound for a trace of the Chern Hessian, or for horizontal derivatives of geometric tensors alone, does
not control this term. Accordingly, \eqref{eq:radial-comparison-assumption} is a sufficient additional assumption for \eqref{eq:radial-B}, not a consequence claimed here of the mixed weighted Ricci bound. The
main theorem only uses \eqref{eq:cutoff-hypotheses}.

\subsection{Weak solutions and regularity}
\label{subsec:weak-solutions}

On a relatively compact smooth domain $\Omega$, the spaces $H^1(\Omega)$ and $H^1_0(\Omega)$ may be defined using any smooth auxiliary Riemannian metric. On compact sets the resulting norms are equivalent to the Finsler
Sobolev norms. We write $H^{-1}(\Omega)$ for the dual of $H^1_0(\Omega)$.

\begin{defi}[Local weak solution]\label{def:weak-solution}
A positive function $u$ is a local weak solution of \eqref{1.1} on $M\times I$ if $u\in L^2_{\mathrm{loc}}(I;H^1_{\mathrm{loc}}(M))$, $u_t\in
L^2_{\mathrm{loc}}(I;H^{-1}_{\mathrm{loc}}(M))$, the reaction $au\log u+bu$ is locally integrable, and
\[
    -\int_I\int_M u\varphi_t\,d\mu\,dt +\int_I\int_M d\varphi(\nabla u)\,d\mu\,dt =\int_I\int_M(au\log u+bu)\varphi\,d\mu\,dt
\]
for every $\varphi\in C_c^\infty(M\times I)$.
\end{defi}

For the gradient estimates, we use solutions with the local interior regularity
\begin{equation}\label{eq:solution-regularity}
\begin{split}
    &u,\ d_xu,\ u_t\ \text{continuous},\\
    &u\in L^2_{\mathrm{loc}}(I;H^2_{\mathrm{loc}}(M)), \qquad u_t\in L^2_{\mathrm{loc}}(I;H^1_{\mathrm{loc}}(M)).
\end{split}
\end{equation}
These are interior assumptions and impose no boundary trace. The solutions constructed in Appendix~\ref{app:existence} have this regularity at positive times, as explained in
Proposition~\ref{prop:interior-regularity}. The H\"older exponent in the local regularity theory depends on the ellipticity and the local coefficient bounds. Smoothness across
$\{du=0\}$ is not required.

Time convolution may be used to justify distributional identities, but it does not preserve the nonlinear equation. If $u_\varepsilon=J_\varepsilon*u$ is a time convolution on
a smaller interval, then
\[
    \partial_tu_\varepsilon =\operatorname{div}_\mu(J_\varepsilon*\nabla u) +J_\varepsilon*(au\log u+bu).
\]
In general $J_\varepsilon*\nabla u\ne\nabla u_\varepsilon$, and the convolved reaction is not $a u_\varepsilon\log u_\varepsilon+b u_\varepsilon$. The proof below differentiates the weak equation by
difference quotients and uses the integrated Bochner inequality, rather than treating time averages as new solutions.

\section{Proof of the local gradient estimate}
\label{sec:gradient-proof}

Throughout this section, all coefficient bounds are taken on $B_p(2R)\times(0,T]$. We write
\[
    f=\log\frac{u}{D},\qquad w=F^2(\nabla f),\qquad L=t\{w+(A+a)f+2(E+b)-2f_t\}.
\]
Since $u>0$, positive homogeneity gives $\nabla f=u^{-1}\nabla u$ and $g_{\nabla f}=g_{\nabla u}$ wherever $du\ne0$. In particular, $\Delta^{\nabla u}f=\Delta f$. The equation for $u$ becomes
\begin{equation}\label{eq:log-equation}
    f_t=\Delta f+w+af+a\log D+b.
\end{equation}
We shall use the following two identities,
\begin{align}
    \partial_t w &=2\langle\nabla f,\nabla^{\nabla u}f_t\rangle_{\nabla u},
\label{eq:time-w}\\
    f_{tt} &=\Delta^{\nabla u}f_t +2\langle\nabla f,\nabla^{\nabla u}f_t\rangle_{\nabla u} +af_t+a_t(f+\log D)+b_t.
\label{eq:time-f}
\end{align}
The second identity follows by differentiating \eqref{eq:log-equation}; the term $a_t f$ is part of the derivative of $af$. The derivative of the Legendre map is the inverse
fundamental tensor away from the zero section, which gives $\partial_t(\Delta f)=\Delta^{\nabla u}f_t$.

The weighted Bochner inequality of Ohta and Sturm \cite{OSSKT2014} gives
\begin{equation}\label{eq:bochner-used}
    \frac12\Delta^{\nabla u}w -d(\Delta f)(\nabla f) \ge -Kw+\frac{(\Delta f)^2}{N}.
\end{equation}
Here the mixed weighted Ricci lower bound is used with equal reference and curvature directions, so it implies $\operatorname{Ric}^N(\nabla f)\ge-KF^2(\nabla f)$. All Laplacians
acting on $a$, $b$, $f_t$, and $L$ below are linearized with reference vector $\nabla u$.

\begin{lemm}\label{Theorem3.0}
For every constant $E$, the quantity $L$ satisfies
\begin{align}\label{eq:differential-L}
    &\Delta^{\nabla u}L-L_t +2\langle\nabla f,\nabla^{\nabla u}L\rangle_{\nabla u} \notag\\
    &\quad\ge -\frac{L}{t}-aL +t\left\{ (A-2K-2-|\log D|)w+\frac{2(\Delta f)^2}{N} +f(\Delta^{\nabla u}a+a_t) \right\} \notag\\
    &\qquad+t\left\{ 2a(E+b)-(A+a)(a\log D+b) +2a_t\log D+2\Delta^{\nabla u}b \right\} \notag\\
    &\qquad-t\left\{ F_{\nabla u}^2(\nabla^{\nabla u}b) +(1+|\log D|) F_{\nabla u}^2(\nabla^{\nabla u}a) \right\}.
\end{align}
\end{lemm}

\begin{proof}
Apply \eqref{eq:bochner-used} to the term $w$ in $L$ and use the product rule for $\Delta^{\nabla u}$. Subtracting the time derivative, then using
\eqref{eq:log-equation}--\eqref{eq:time-f}, gives
\begin{align*}
    &\Delta^{\nabla u}L-L_t +2\langle\nabla f,\nabla^{\nabla u}L\rangle_{\nabla u} \\
    &\quad\ge -\frac{L}{t}-aL +t\left\{ (A-2K)w+\frac{2(\Delta f)^2}{N} +f(\Delta^{\nabla u}a+a_t) \right\} \\
    &\qquad+t\left\{ 2a(E+b)-(A+a)(a\log D+b) +2a_t\log D+2\Delta^{\nabla u}b \right\} \\
    &\qquad+2t\langle\nabla f,\nabla^{\nabla u}b\rangle_{\nabla u} +2t(1-\log D) \langle\nabla f,\nabla^{\nabla u}a\rangle_{\nabla u}.
\end{align*}
For completeness, the term $-aL$ follows from the identity
\[
    2af_t-a(A+a)f =aw+2a(E+b)-\frac{aL}{t}.
\]
The derivatives of $a$ in the Bochner term and the product rule combine as
\[
    -2(f+\log D)d a(\nabla f)+2(f+1)d a(\nabla f) =2(1-\log D)d a(\nabla f).
\]
Finally, the Cauchy--Schwarz inequality for $g_{\nabla u}$ gives
\begin{align*}
    2\langle\nabla f,\nabla^{\nabla u}b\rangle_{\nabla u} &\ge-w-F_{\nabla u}^2(\nabla^{\nabla u}b),\\
    2(1-\log D) \langle\nabla f,\nabla^{\nabla u}a\rangle_{\nabla u} &\ge-(1+|\log D|) \{w+F_{\nabla u}^2(\nabla^{\nabla u}a)\}.
\end{align*}
This proves \eqref{eq:differential-L}.
\end{proof}

We choose $E\ge0$ so that
\begin{equation}\label{eq:E-positive}
    E+b\ge0,\qquad E-a\log D\ge0,
\end{equation}
and
\begin{align}\label{eq:E-condition}
    &\frac{2(E-a\log D)^2}{N} +2a(E+b)-(A+a)(a\log D+b) +2a_t\log D+2\Delta^{\nabla u}b \notag\\
    &\qquad -F_{\nabla u}^2(\nabla^{\nabla u}b) -(1+|\log D|) F_{\nabla u}^2(\nabla^{\nabla u}a) \ge0.
\end{align}
Such a choice is possible under the stated coefficient bounds. Indeed, the first term grows quadratically in $E$, the possibly negative term $2a(E+b)$ grows at most linearly, and all remaining
terms have a uniform lower bound. This also explains the dependence of $E$ on $A$, $D$, $N$, and the coefficient bounds.

The definition of $L$ and \eqref{eq:log-equation} give
\begin{equation}\label{eq:Delta-f-L}
    \Delta f =-\frac{L}{2t}-\frac{w}{2} +\frac{A-a}{2}f+E-a\log D.
\end{equation}
In particular, the terms containing $b$ cancel in this identity. On the set where $L>0$, put $h=w/L\ge0$. Squaring \eqref{eq:Delta-f-L} and discarding the nonnegative term $(A-a)^2f^2/4$
gives
\begin{align*}
    (\Delta f)^2 &\ge\frac{(1+ht)^2L^2}{4t^2}+(E-a\log D)^2 -\frac{(1+ht)(A-a)Lf}{2t}\\
    &\qquad-\frac{(1+ht)(E-a\log D)L}{t} +(A-a)(E-a\log D)f.
\end{align*}
Consequently, \eqref{eq:differential-L} and \eqref{eq:E-condition} imply
\begin{align}\label{eq:L-reduced}
    &\Delta^{\nabla u}L-L_t +2\langle\nabla f,\nabla^{\nabla u}L\rangle_{\nabla u} \notag\\
    &\quad\ge -\frac{L}{t}-aL +t\left\{ (A-2K-2-|\log D|)hL +\frac{(1+ht)^2L^2}{2Nt^2} \right\} \notag\\
    &\qquad-\frac{2(1+ht)(E-a\log D)L}{N} \notag\\
    &\qquad+tf\left\{ -\frac{(A-a)(1+ht)L}{Nt} +\Delta^{\nabla u}a+a_t +\frac{2(A-a)(E-a\log D)}{N} \right\}.
\end{align}

\begin{proof}[Proof of Theorem~\ref{MainTheorem}]
Let $\phi$ be the cutoff in the assumptions and put $H=\phi L$. We first give the computation for regular points; the weak formulation and the treatment of the critical
set are explained below. If the maximum of $H$ on $B_p(2R)\times[0,T]$ is nonpositive, the conclusion follows from $E+b\ge0$. Otherwise, at a positive maximum $(x_0,t_0)$, one has
\[
    dH=0,\qquad H_t\ge0,\qquad\Delta^{\nabla u}H\le0.
\]
The product rule and the cutoff assumptions imply
\begin{equation}\label{eq:cutoff-at-max}
    \phi(\Delta^{\nabla u}L-L_t) \le L\left\{ \frac{2F_{\nabla u}^2(\nabla^{\nabla u}\phi)}{\phi} -\Delta^{\nabla u}\phi \right\} \le BL.
\end{equation}
All quantities in the next computation are evaluated at $(x_0,t_0)$, and we omit the subscripts.

Suppose first that the expression in braces multiplying $tf$ in \eqref{eq:L-reduced} is nonnegative. Since $A>[a]^+$, it follows that
\begin{align}\label{eq:first-case}
    H\le L &\le\frac{Nt}{(A-a)(1+ht)} \left\{\Delta^{\nabla u}a+a_t +\frac{2(A-a)(E-a\log D)}{N}\right\} \notag\\
    &\le t\left\{ \frac{N[\Delta^{\nabla u}a+a_t]^+}{A-[a]^+} +2[E-a\log D]^+ \right\}.
\end{align}
If that expression is negative, then its product with $f\le0$ is nonnegative and can be discarded. Moreover,
\[
    -2\phi\langle\nabla f,\nabla^{\nabla u}L\rangle_{\nabla u} =2L\langle\nabla f,\nabla^{\nabla u}\phi\rangle_{\nabla u} \ge-\frac{2\sqrt{\alpha}\,C_1}{R} \phi^{1/2}L^{3/2}h^{1/2}.
\]
Using \eqref{eq:cutoff-at-max}, multiplying by $t/L$, and using $0\le\phi\le1$, we obtain
\begin{align*}
    \frac{(1+ht)^2H}{2N} &\le 1+Bt+[a]^+t +\frac{2t(1+ht)[E-a\log D]^+}{N}\\
    &\qquad+t^2h(2K+2+|\log D|-A)^+ +\frac{2\sqrt{\alpha}\,C_1t}{R}\sqrt{Hh}.
\end{align*}
Young's inequality gives
\[
    \frac{2\sqrt{\alpha}\,C_1t}{R}\sqrt{Hh} \le\frac{(1+ht)^2H}{4N} +\frac{4N\alpha C_1^2ht^2}{R^2(1+ht)^2}.
\]
For $y\ge0$, we have $y/(1+y)^2\le1/4$ and $4y/(1+y)^4\le1$. Therefore,
\begin{align}\label{eq:second-case}
    H\le4N\bigg\{1+t\bigg( [a]^++B+\frac{N\alpha C_1^2}{R^2} +\frac{(2K+2+|\log D|-A)^+}{4} \notag\\
    {}+\frac{2[E-a\log D]^+}{N} \bigg)\bigg\}.
\end{align}
The right-hand sides of \eqref{eq:first-case} and \eqref{eq:second-case} are nondecreasing in $t$. Since $t_0\le T$ and $\phi=1$ on $B_p(R)$, they yield the slightly weaker common bound
\begin{align*}
    &w+(A+a)f-2f_t\\
    &\quad\le4N\left\{ \frac1T+[a]^++B+\frac{N\alpha C_1^2}{R^2} \right.\\
    &\qquad\left. +\frac{(2K+2+|\log D|-A)^+}{4} +\frac{2[E-a\log D]^+}{N} +\frac{[\Delta^{\nabla u}a+a_t]^+}{4(A-[a]^+)} \right\}
\end{align*}
on $B_p(R)\times\{T\}$. Here we also used $w+(A+a)f-2f_t=L/T-2(E+b)\le L/T$.

No initial derivative bound is needed. For $\varepsilon>0$, replace $t$ by $t-\varepsilon$ in $L$ and apply the argument on $[\varepsilon,T]$. Then $L=0$ at the
initial time. The resulting bound has $1/(T-\varepsilon)$ in place of $1/T$; letting $\varepsilon\downarrow0$ proves the stated estimate.

We finally justify the argument at critical points of $u$. Extend $g_{\nabla u}$ measurably and uniformly elliptically across $\{du=0\}$, as in the definition of the
linearized operator. The regularity assumptions and \eqref{eq:log-equation} give $\Delta f\in L^2_{\mathrm{loc}}((0,T];H^1_{\mathrm{loc}})$. The integrated Bochner
inequality \cite{OSSKT2014} therefore gives \eqref{eq:bochner-used} on the whole cylinder. The chain and product rules hold weakly. In particular, the locally Lipschitz Legendre map can be differentiated almost everywhere; on the
set $df=0$, its time derivative and $df_t$ vanish almost everywhere. Thus \eqref{eq:time-w} and \eqref{eq:time-f} also hold weakly.

It follows that \eqref{eq:L-reduced} is a weak differential inequality on the whole cylinder. To apply the weak maximum principle, work on a positive superlevel
set of $H$. On $\{\phi>0\}$ the product rule gives
\begin{align*}
    &(\Delta^{\nabla u}-\partial_t)H +2dH(\nabla f) -\frac{2}{\phi}g_{\nabla u}^*(d\phi,dH)\\
    &\quad=\phi\bigl\{ (\Delta^{\nabla u}-\partial_t)L+2dL(\nabla f) \bigr\}\\
    &\qquad+L\left\{ \Delta^{\nabla u}\phi -\frac{2F_{\nabla u}^2(\nabla^{\nabla u}\phi)}{\phi} \right\} +2L\,d\phi(\nabla f).
\end{align*}
The cutoff assumption bounds the expression in the second pair of braces from below by $-B$. If $H$ exceeds the right-hand side of \eqref{eq:first-case} with $t$ replaced by $T$, the expression multiplying $tf$ in
\eqref{eq:L-reduced} is nonpositive throughout that set. Hence the second-case inequality applies there. Testing the resulting inequality with a positive part of $H$
above a proposed bound, with a small strict time barrier, gives \eqref{eq:second-case}. The cutoff terms are controlled by the assumed weak cutoff inequalities. This is the
weak maximum principle for a uniformly elliptic divergence-form operator on a compact cylinder; it does not require the maximum to lie in $\{du\ne0\}$. This proves the estimate with the
regularity stated in the theorem.
\end{proof}

\begin{coro}\label{cor:linear-gradient}
Assume the curvature and cutoff hypotheses of Theorem~\ref{MainTheorem}. Let $u>0$ have the interior regularity \eqref{eq:solution-regularity} and solve
\[
    u_t=\Delta u+bu.
\]
Assume that $b$ is $C^2$ in space and $C^1$ in time and that $\Delta^{\nabla u}b$ has a locally integrable representative, as in Theorem~\ref{MainTheorem}. No upper bound on $u$ is needed. Suppose $E\ge0$
satisfies
\begin{equation}\label{eq:E-linear}
    E+b\ge0,\qquad \frac{2E^2}{N}+2\Delta^{\nabla u}b -F_{\nabla u}^2(\nabla^{\nabla u}b)\ge0.
\end{equation}
Then
\begin{equation}\label{eq:linear-gradient}
    \frac{F^2(\nabla u)}{u^2}-2\frac{u_t}{u} \le4N\left\{ \frac1t+B+\frac{N\alpha C_1^2}{R^2} +\frac{2K+1}{4}+\frac{2E}{N} \right\}
\end{equation}
on $B_p(R)\times(0,T]$.
\end{coro}

\begin{proof}
Use $f=\log u$ and $L=t\{w+2(E+b)-2f_t\}$. The terms involving $f$ and derivatives of $a$ are absent. The only coefficient-gradient term is $2\langle\nabla
f,\nabla^{\nabla u}b\rangle_{\nabla u}$, so the coefficient of $w$ after Young's inequality is $-2K-1$. Furthermore,
\[
    \Delta f=-\frac{L}{2t}-\frac{w}{2}+E.
\]
Condition \eqref{eq:E-linear} leaves precisely the inequality used in the second case above, with $a=0$ and curvature term $2K+1$. The same cutoff argument proves
\eqref{eq:linear-gradient}. This proof uses neither $f\le0$ nor the condition $A>[a]^+$ from the nonlinear theorem.
\end{proof}
\section{Applications}\label{sec:applications}

We first integrate the gradient estimate to obtain a local Harnack inequality for the equation with $a=0$. We then prove an upper bound for positive solutions of a stationary logarithmic equation. The latter argument does
not assume that the solution is bounded.

\begin{theo}[Local Harnack inequality]\label{thm:harnack}
Assume the hypotheses of Corollary~\ref{cor:linear-gradient}, and let $u>0$ solve
\[
    u_t=\Delta u+b(x,t)u
\]
on $B_p(2R)\times(0,T]$. Choose $E\geq0$ so that
\[
    E+b\geq0, \qquad \frac{2E^2}{N}+2\Delta^{\nabla u}b -F_{\nabla u}^{2}(\nabla^{\nabla u}b)\geq0
\]
on this cylinder. For $x_1,x_2\in B_p(R)$ and $0<t_1<t_2\leq T$,
\begin{align}\label{eq:harnack}
    u(x_1,t_1) &\leq u(x_2,t_2) \left(\frac{t_2}{t_1}\right)^{2N} \exp\left\{ 2N(t_2-t_1) \left(B+\frac{N\alpha C_1^2}{R^2} +\frac{2K+1}{4}+\frac{2E}{N}\right) \right.\notag\\
    &\hspace{37mm}\left. +\inf_{\gamma} \frac{1}{2(t_2-t_1)} \int_0^1 F^2(\dot\gamma(s))\,ds \right\}.
\end{align}
The infimum is taken over piecewise $C^1$ curves in $B_p(R)$ satisfying $\gamma(0)=x_2$ and $\gamma(1)=x_1$.
\end{theo}

\begin{proof}
Put $f=\log u$, fix an admissible curve $\gamma$, and write $t(s)=(1-s)t_2+st_1$. Along this space-time curve,
\[
    \log\frac{u(x_1,t_1)}{u(x_2,t_2)} =\int_0^1 \bigl[df(\dot\gamma)-(t_2-t_1)f_t\bigr]\,ds.
\]
The duality inequality gives $df(\dot\gamma)\leq F(\nabla f)F(\dot\gamma)$. By Corollary~\ref{cor:linear-gradient} and the inequality
\[
    F(\nabla f)F(\dot\gamma) \leq\frac{t_2-t_1}{2}F^2(\nabla f) +\frac{F^2(\dot\gamma)}{2(t_2-t_1)},
\]
we obtain
\begin{align*}
    \log\frac{u(x_1,t_1)}{u(x_2,t_2)} &\leq\int_0^1 \frac{F^2(\dot\gamma(s))}{2(t_2-t_1)}\,ds +2N\int_0^1\frac{t_2-t_1}{t(s)}\,ds\\
    &\quad+2N(t_2-t_1) \left(B+\frac{N\alpha C_1^2}{R^2} +\frac{2K+1}{4}+\frac{2E}{N}\right).
\end{align*}
Here
\[
    \int_0^1\frac{t_2-t_1}{(1-s)t_2+st_1}\,ds =\log\frac{t_2}{t_1}.
\]
Taking the infimum over $\gamma$ and exponentiating proves \eqref{eq:harnack}.
\end{proof}

The direction of the curves in \eqref{eq:harnack} matters when $F$ is not reversible. The argument uses curves from the point at the later time to the point at the earlier time. This choice
permits the direct use of $df(\dot\gamma)\leq F^*(df)F(\dot\gamma)$.

We next consider the stationary equation. The cutoff assumption is stated explicitly because the limit as $R\to\infty$ requires control of the divergence-form linearization
along the solution.

\begin{theo}[An upper bound for stationary solutions]
\label{thm:stationary}
Let $(M,F,\mu)$ be a forward complete noncompact Finsler manifold with finite misalignment and $\mathrm{Ric}^N\geq0$, where $n<N<\infty$. Let $u>0$ be a
local weak solution with $u\in H^2_{\mathrm{loc}}(M)\cap C^1(M)$ of
\begin{equation}\label{eq:stationary}
    \Delta u+2u\log u+V(x)u=0,
\end{equation}
where $V\in C^2(M)$. Assume that there are compactly supported cutoffs $0\leq\phi_R\leq1$, equal to one on $B_p(R)$, for which
\begin{equation}\label{eq:stationary-cutoffs}
    \frac{F_{\nabla u}^{2}(\nabla^{\nabla u}\phi_R)}{\phi_R} \leq\frac{\alpha C_1^2}{R^2}, \qquad \frac{2F_{\nabla u}^{2}(\nabla^{\nabla u}\phi_R)}{\phi_R} -\Delta^{\nabla
    u}\phi_R\leq B_R,
\end{equation}
on $\{\phi_R>0\}$, with $B_R\to0$ as $R\to\infty$. Suppose $\Delta^{\nabla u}V$ has a locally integrable representative and that $V$, $F_{\nabla u}(\nabla^{\nabla u}V)$, and the negative part of $\Delta^{\nabla u}V$
are bounded. If $E\geq0$ satisfies
\begin{equation}\label{eq:stationary-E}
    E+V\geq0, \qquad \frac{2E^2}{N}+2\Delta^{\nabla u}V -F_{\nabla u}^{2}(\nabla^{\nabla u}V)\geq0,
\end{equation}
then
\begin{equation}\label{eq:stationary-bound}
    u\leq\exp(2N+2E) \qquad\text{on }M.
\end{equation}
For example, one may take
\begin{align*}
    E={}&\sup_M|V| +\sqrt{N\sup_M\max\{-\Delta^{\nabla u}V,0\}}\\
    &+\sqrt{\frac N2} \sup_M F_{\nabla u}(\nabla^{\nabla u}V).
\end{align*}
\end{theo}

\begin{proof}
Set
\[
    f=\log u, \qquad w=F^2(\nabla f), \qquad Q=w+4f+2(V+E).
\]
No sign condition on $f$ is imposed. Equation~\eqref{eq:stationary} gives
\[
    \Delta f=-w-2f-V=-\frac{Q+w}{2}+E.
\]
The Bochner inequality under $\mathrm{Ric}^N\geq0$ yields
\begin{align*}
    \Delta^{\nabla f}Q+2\,dQ(\nabla f) &\geq\frac{2}{N}(\Delta f)^2 -2Q+2w+4E +2\,dV(\nabla f)+2\Delta^{\nabla f}V.
\end{align*}
Indeed, the terms involving $dw(\nabla f)$ cancel after adding $2\,dQ(\nabla f)$. Since
\[
    2\,dV(\nabla f) \geq-w-F_{\nabla f}^{2}(\nabla^{\nabla f}V),
\]
conditions~\eqref{eq:stationary-E} imply
\begin{equation}\label{eq:stationary-Q}
    \Delta^{\nabla f}Q+2\,dQ(\nabla f) \geq\frac{(Q+w)^2}{2N} -\frac{2E}{N}(Q+w)-2Q+w.
\end{equation}
We use the cutoff $\phi_R$ and consider a positive maximum of $\phi_RQ$. If this maximum is not positive, the desired local upper bound is immediate.
Otherwise, put $P=\phi_RQ$ and $v=w/Q$ at a positive maximum point. Thus $v\geq0$ and $dQ=-Q\,d\phi_R/\phi_R$ there. The maximum principle, \eqref{eq:stationary-cutoffs}, and \eqref{eq:stationary-Q} give
\[
    \frac{(1+v)^2}{2N}P \leq2+B_R+\frac{2E}{N}(1+v) +2\sqrt{\frac{\alpha C_1^2}{R^2}\,vP}.
\]
This calculation can equivalently be performed with positive-part test functions for the weak Bochner inequality, as in the proof of
Theorem~\ref{MainTheorem}. By Young's inequality and $4v\leq(1+v)^2$,
\[
    2\sqrt{\frac{\alpha C_1^2}{R^2}\,vP} \leq\frac{(1+v)^2}{4N}P +\frac{N\alpha C_1^2}{R^2}.
\]
It follows that
\[
    \phi_RQ \leq4N\left( 2+B_R+\frac{N\alpha C_1^2}{R^2}+\frac{2E}{N} \right)
\]
everywhere. On $B_p(R)$, where $\phi_R=1$, let $R\to\infty$ to obtain $Q\leq8N+8E$. Finally, $w\geq0$ and $E+V\geq0$ show that $4\log u\leq Q$, which proves~\eqref{eq:stationary-bound}.
\end{proof}
\appendix
\section{Existence of positive solutions}
\label{app:existence}

We give a variational construction for the equation studied in this paper. The metric and the measure are fixed, while the coefficients may depend on both space and time. The construction uses the
convex Finsler energy on the Hilbert space of square-integrable functions, as in the heat-flow construction of Ohta and Sturm \cite{OhtaSturm2009}. The logarithmic
reaction requires an additional truncation and comparison argument. On a noncompact manifold, we work first on relatively compact domains and then pass to an exhaustion. This
avoids an unnecessary global square-integrability assumption on positive bounded solutions.

\begin{theo}[Existence for positive bounded initial data]
\label{thm:existence}
Let $(M,F,\mu)$ be a smooth forward complete Finsler metric measure manifold without boundary. Suppose that $a$ and $b$ have the regularity assumed in the main theorem
and satisfy
\[
    \sup_{M\times(0,T)}|a|+\sup_{M\times(0,T)}|b|<\infty \qquad\text{for every }T>0.
\]
Let $u_0$ be a measurable function satisfying
\[
    0<c\leq u_0\leq C<\infty \qquad\text{almost everywhere on }M.
\]
Then there exists a positive weak solution of
\begin{equation}\label{app:cauchy}
    u_t=\Delta u+a(x,t)u\log u+b(x,t)u, \qquad u(\cdot,0)=u_0,
\end{equation}
on $M\times(0,\infty)$. For every $T>0$ and every relatively compact smooth domain $\Omega\subset M$, it satisfies
\[
\begin{split}
    u&\in C([0,T];L^2(\Omega))\cap L^2(0,T;H^1(\Omega)),\\
    u_t&\in L^2(0,T;H^{-1}(\Omega)).
\end{split}
\]
Here $H^{-1}(\Omega)$ is the dual of $H^1_0(\Omega)$. More explicitly, put
\[
    A=\sup_{M\times(0,T)}|a|,\qquad B=\sup_{M\times(0,T)}|b|,\qquad z(0)=\max\{0,-\log c,\log C\},
\]
and let $z'=Az+B$. Then
\begin{equation}\label{app:positive-bounds}
    e^{-z(t)}\leq u(x,t)\leq e^{z(t)} \qquad\text{for almost every }(x,t)\in M\times(0,T).
\end{equation}
No curvature bound is needed for this existence result.
\end{theo}

\begin{proof}
Fix $T>0$. We use an auxiliary smooth Riemannian metric to define Sobolev norms on relatively compact domains. These norms are equivalent to the corresponding Finsler Sobolev norms on each such
domain.

\medskip \noindent\emph{The convex energy.} For a relatively compact smooth domain $\Omega$, set
\[
    \mathcal E(v)=\frac12\int_\Omega F^*(dv)^2\,d\mu, \qquad v\in H^1_0(\Omega).
\]
Extend $\mathcal E$ by $+\infty$ outside $H^1_0(\Omega)$ in $L^2(\Omega,\mu)$. Strong convexity of $F$ implies that $\mathcal E$ is convex and lower
semicontinuous. If
\[
    J^*(x,\xi)=\partial_\xi\left(\frac12F^*(x,\xi)^2\right), \qquad J^*(x,0)=0,
\]
then $J^*(du)=\nabla u$, and differentiation in the direction $\phi\in H^1_0(\Omega)$ gives
\begin{equation}\label{app:variation}
    \left.\frac{d}{ds}\mathcal E(v+s\phi)\right|_{s=0} =\int_\Omega d\phi(\nabla v)\,d\mu.
\end{equation}
Consequently the $L^2$ subgradient of $\mathcal E$ is $-\Delta$ on its domain. The Hilbert structure used here is
\[
    (v,w)_{L^2}=\int_\Omega vw\,d\mu.
\]
It does not require the Finsler norm on $T_xM$ to come from an inner product. In general, $\mathcal E$ is not quadratic and its subgradient is nonlinear.

For later use, smoothness and strong convexity give positive constants $\lambda_\Omega$ and $\Lambda_\Omega$ such that
\begin{equation}\label{app:monotonicity}
\begin{split}
    (\xi-\eta)\bigl(J^*(x,\xi)-J^*(x,\eta)\bigr) &\geq\lambda_\Omega|\xi-\eta|^2,\\
    |J^*(x,\xi)-J^*(x,\eta)| &\leq\Lambda_\Omega|\xi-\eta|
\end{split}
\end{equation}
for $x\in\overline\Omega$ and all covectors $\xi,\eta\in T_x^*M$. Indeed, the fiber Hessian of $F^{*2}/2$ is positive definite and homogeneous of degree zero away from the zero section. Compactness of
the unit cotangent sphere over $\overline\Omega$ gives uniform bounds for that Hessian. Integration along the line segment from $\eta$ to $\xi$ proves \eqref{app:monotonicity}, also when the segment passes through zero.

\medskip \noindent\emph{A bounded-domain problem with a truncated reaction.} Write $\ell=e^{-z(T)}$ and $L=e^{z(T)}$, and define
\[
    P(s)=\min\{L,\max\{\ell,s\}\},\qquad G(x,t,s)=a(x,t)P(s)\log P(s)+b(x,t)P(s).
\]
The map $s\mapsto G(x,t,s)$ is bounded and globally Lipschitz, uniformly in $(x,t)\in M\times(0,T)$. A Lipschitz constant is $A(1+z(T))+B$. We first solve
\begin{equation}\label{app:truncated}
    u_t=\Delta u+G(x,t,u),\qquad u|_{\partial\Omega}=1,\qquad u(\cdot,0)=u_0|_\Omega.
\end{equation}
Put $v=u-1$, so that $v$ has zero boundary values and $\nabla u=\nabla v$.

Here is a direct construction. Divide $[0,T]$ into intervals of length $h=T/m$. Choose $v_0^h\in H^1_0(\Omega)$ such that
\[
    v_0^h\longrightarrow u_0-1\quad\text{in }L^2(\Omega), \qquad h\|v_0^h\|_{H^1}^2\leq C_\Omega.
\]
Such a choice follows from density, by letting the approximation index tend to infinity sufficiently slowly as $h$ tends to zero. Having chosen $v_{i-1}$,
define
\[
    G_i(x,s)=\frac1h\int_{(i-1)h}^{ih}G(x,t,s+1)\,dt
\]
and let $v_i$ minimize
\begin{equation}\label{app:minimization}
    \mathcal E(v)+\frac{1}{2h}\|v-v_{i-1}\|_{L^2(\Omega)}^2 -\int_\Omega G_i(x,v_{i-1})v\,d\mu
\end{equation}
over $H^1_0(\Omega)$. Coercivity and weak lower semicontinuity give a minimizer, and the quadratic term makes it unique. Its Euler equation is
\begin{equation}\label{app:euler}
    \int_\Omega\frac{v_i-v_{i-1}}{h}\phi\,d\mu +\int_\Omega d\phi(\nabla v_i)\,d\mu =\int_\Omega G_i(x,v_{i-1})\phi\,d\mu.
\end{equation}

Taking $\phi=v_i$ and using
\[
    2(v_i-v_{i-1},v_i)_{L^2} =\|v_i\|_{L^2}^2-\|v_{i-1}\|_{L^2}^2 +\|v_i-v_{i-1}\|_{L^2}^2
\]
gives, by Young's inequality and the discrete Gronwall inequality,
\begin{equation}\label{app:discrete-bounds}
    \max_i\|v_i\|_{L^2}^2 +\sum_i\|v_i-v_{i-1}\|_{L^2}^2 +h\sum_i\|v_i\|_{H^1}^2\leq C_{\Omega,T}.
\end{equation}
The constant is independent of $h$. In obtaining this estimate we use $dv_i(\nabla v_i)=F^*(dv_i)^2$ and the boundedness of $G$. Equation \eqref{app:euler} and \eqref{app:monotonicity} also give
\[
    h\sum_i\left\|\frac{v_i-v_{i-1}}h\right\|_{H^{-1}(\Omega)}^2 \leq C_{\Omega,T}.
\]
Let $v_h=v_i$ on $((i-1)h,ih]$, and let $\widetilde v_h$ be the piecewise affine interpolation. The preceding estimates bound $\widetilde v_h$ in $L^2(0,T;H^1_0(\Omega))$ and its time derivative in
$L^2(0,T;H^{-1}(\Omega))$. The compact embedding $H^1_0(\Omega)\hookrightarrow L^2(\Omega)$ therefore gives strong compactness in $L^2(\Omega\times(0,T))$ by the Aubin--Lions
compactness theorem \cite{Simon1987}. Furthermore,
\[
    \|v_h-\widetilde v_h\|_{L^2(\Omega\times(0,T))}^2 \leq h\sum_i\|v_i-v_{i-1}\|_{L^2}^2\longrightarrow0.
\]
The interpolation equal to $v_{i-1}$ on the same intervals has the same limit. After taking a subsequence, $v_h$ converges strongly to $v$ in $L^2$, $dv_h$
converges weakly to $dv$, and $J^*(dv_h)$ converges weakly to a vector field $q$. The time averages of $a$ and $b$ converge to the corresponding coefficients in $L^2(\Omega\times(0,T))$. This fact, together with the
boundedness and Lipschitz continuity of $P$ and $P\log P$, shows that the averaged reaction converges strongly to $G(x,t,v+1)$. The affine interpolations are equicontinuous in $H^{-1}(\Omega)$ and
uniformly bounded in $L^2(\Omega)$. Testing against $H^1_0(\Omega)$ and then using density identifies their weak $L^2$ traces at every time. In particular, $v_h(T)\rightharpoonup v(T)$ in
$L^2(\Omega)$.

We verify the nonlinear diffusion term. Passing to the limit in \eqref{app:euler} first gives
\[
    v_t=\operatorname{div}_\mu q+G(x,t,v+1), \qquad v(0)=u_0-1.
\]
The usual Hilbert-space energy identity applies to $v\in L^2(0,T;H^1_0(\Omega))$ with $v_t\in L^2(0,T;H^{-1}(\Omega))$; in particular, $v\in C([0,T];L^2(\Omega))$. It gives
\[
    \frac12\|v(T)\|_{L^2}^2-\frac12\|v(0)\|_{L^2}^2 =\int_0^T\langle v_t,v\rangle_{H^{-1},H^1_0}\,dt.
\]
Summing \eqref{app:euler} with $\phi=v_i$ and using weak lower semicontinuity of the terminal $L^2$ norm gives
\[
    \limsup_{h\to0}\int_0^T\int_\Omega dv_h\bigl(J^*(dv_h)\bigr)\,d\mu\,dt \leq\int_0^T\int_\Omega dv(q)\,d\mu\,dt.
\]
For every square-integrable covector field $\xi$, monotonicity gives
\[
    \int_0^T\int_\Omega (dv_h-\xi)\bigl(J^*(dv_h)-J^*(\xi)\bigr)\,d\mu\,dt\geq0.
\]
The preceding upper bound allows passage to the limit in this inequality. Choosing $\xi=dv+s\eta$ for an arbitrary square-integrable covector field $\eta$, dividing by $s$ with each sign, and letting $s$ tend to zero yields
$q=J^*(dv)$. Thus $u=v+1$ solves \eqref{app:truncated}.

\medskip \noindent\emph{Positive upper and lower bounds.} The spatially constant functions $e^{-z(t)}$ and $e^{z(t)}$ satisfy
\[
\begin{split}
    \frac{d}{dt}e^{-z(t)} &=(-Az-B)e^{-z(t)} \leq a(x,t)e^{-z(t)}\log(e^{-z(t)})+b(x,t)e^{-z(t)},\\
    \frac{d}{dt}e^{z(t)} &=(Az+B)e^{z(t)} \geq a(x,t)e^{z(t)}\log(e^{z(t)})+b(x,t)e^{z(t)}.
\end{split}
\]
Their values lie in $[\ell,L]$, so these inequalities also hold for the truncated reaction. They bound the initial data and the lateral boundary value $1$. Testing
the difference between \eqref{app:truncated} and the upper inequality with $(u-e^{z(t)})_+$, and using monotonicity and the Lipschitz bound for $G$, gives
\[
    \frac12\frac{d}{dt}\|(u-e^{z(t)})_+\|_{L^2(\Omega)}^2 \leq \bigl(A(1+z(T))+B\bigr) \|(u-e^{z(t)})_+\|_{L^2(\Omega)}^2.
\]
The initial value of this norm is zero. Gronwall's inequality gives $u\leq e^{z(t)}$. Testing with $(e^{-z(t)}-u)_+$ proves the lower bound. Therefore $P(u)=u$, and the truncation disappears from the equation. The
same difference argument proves uniqueness on $\Omega$ among solutions with the same boundary and initial data that stay in a fixed positive bounded interval on each finite time
interval.

\medskip \noindent\emph{Passage to the whole manifold.} Choose a smooth exhaustion $\Omega_1\Subset\Omega_2\Subset\cdots\Subset M$ and construct the
preceding solution $u_j$ on each $\Omega_j$, with boundary value $1$ and initial data $u_0|_{\Omega_j}$. The bounds \eqref{app:positive-bounds} are independent of $j$. For a
compactly supported smooth spatial cutoff $\eta$, testing the equation by $\eta^2(u_j-1)$ and applying Young's inequality gives uniform local bounds for
\[
    \sup_{0\leq t\leq T}\int_M\eta^2|u_j-1|^2\,d\mu +\int_0^T\int_M\eta^2|du_j|^2\,d\mu\,dt.
\]
The equation also bounds $(u_j)_t$ locally in $L^2(0,T;H^{-1})$. Compactness and a diagonal subsequence give strong convergence $u_j\to u$ in $L^2$ on every compact spatial domain times $(0,T)$.

To pass through the nonlinear gradient, fix such a cutoff $\eta$ and take $j,k$ large enough that its support lies in both domains. Subtract the equations and test by $\eta^2(u_j-u_k)$. Using both
inequalities in \eqref{app:monotonicity} on the support of $\eta$ and the Lipschitz bound for the reaction gives
\[
\begin{split}
    &\sup_{0\leq t\leq T}\int_M\eta^2|u_j-u_k|^2\,d\mu +\int_0^T\int_M\eta^2|du_j-du_k|^2\,d\mu\,dt\\
    &\qquad\leq C_{\eta,T}\int_0^T\int_{\operatorname{supp}\eta} |u_j-u_k|^2\,d\mu\,dt.
\end{split}
\]
There is no initial term because both solutions start from $u_0$. The right-hand side tends to zero. Thus the gradients converge strongly locally, and so do
$J^*(du_j)=\nabla u_j$. We may pass to the limit in the weak equation. In particular, for every $\phi\in C_c^\infty(M\times[0,T))$,
\[
\begin{split}
    -\int_0^T\int_M u\phi_t\,d\mu\,dt +\int_0^T\int_M d\phi(\nabla u)\,d\mu\,dt &=\int_0^T\int_M(au\log u+bu)\phi\,d\mu\,dt\\
    &\quad+\int_Mu_0\phi(\cdot,0)\,d\mu.
\end{split}
\]
The same local difference estimate gives continuity into $L^2$ on compact spatial domains, including at $t=0$. The positive bounds pass to the limit as well.

For completeness, the scalar function used above is
\[
    z(t)=
\begin{cases}
    z(0)e^{At}+\dfrac{B}{A}(e^{At}-1),& A>0,\\[4pt]
    z(0)+Bt,& A=0.
\end{cases}
\]
It stays finite on every finite time interval. Uniqueness on each bounded domain makes the constructions on different time intervals consistent. A further diagonal subsequence over $T=1,2,\ldots$
therefore produces a solution for all positive times.
\end{proof}

\begin{rem}\label{rem:existence-scope}
The theorem gives a positive bounded solution on every finite time interval. The upper bound may grow with $T$. A bound valid for all positive times follows, for example, if there are constants $0<c\leq C$ such that the initial
data lie in $[c,C]$ and
\begin{equation}\label{app:invariant-interval}
    a(x,t)\log c+b(x,t)\geq0,\qquad a(x,t)\log C+b(x,t)\leq0
\end{equation}
on $M\times(0,\infty)$. In that case the construction can be chosen to give $c\leq u\leq C$ for all times. One uses a fixed lateral boundary value in $[c,C]$ in
place of $1$ and applies the same comparison argument with the constants $c$ and $C$. Without \eqref{app:invariant-interval}, the bound $D=e^{z(T)}$ is available on a prescribed finite
time interval. The gradient estimates then apply provided their remaining curvature, coefficient, and cutoff hypotheses hold.

If $\mu(M)=\infty$, the lower bound in \eqref{app:positive-bounds} prevents $u(\cdot,t)$ from belonging to $L^2(M,\mu)$. The solution is therefore understood in the local energy sense
throughout $M$, as expressed by the displayed weak formulation. A global $L^2$ energy class requires additional integrability assumptions on the initial data and should not be imposed on the positive class constructed
here.
\end{rem}

\begin{prop}[Interior regularity]\label{prop:interior-regularity}
Under the assumptions of Theorem~\ref{thm:existence}, the solution constructed there has a representative for which $u$, $d_xu$, and $u_t$ are locally H\"older continuous on $M\times(0,\infty)$, with some positive exponent
on each compact interior cylinder. Moreover,
\[
    u\in L^2_{\mathrm{loc}}((0,\infty);H^2_{\mathrm{loc}}(M)), \qquad u_t\in L^2_{\mathrm{loc}}((0,\infty);H^1_{\mathrm{loc}}(M)).
\]
If $a$ and $b$ are smooth, then $u$ is smooth on the open set $\{(x,t)\in M\times(0,\infty)\mid d_xu(x,t)\neq0\}$.
\end{prop}

\begin{proof}
We explain how the interior argument for the Finsler heat equation extends to the present reaction term. Fix a compact interior cylinder and a slightly
larger coordinate cylinder. The bounds \eqref{app:positive-bounds} place $u$ in a fixed positive bounded interval there. Consequently, for
\[
    r(x,t,s)=a(x,t)s\log s+b(x,t)s,
\]
the functions $r$, $r_s$, $r_t$, and $d_xr$ are bounded on the relevant range of $(x,t,s)$. In the coordinates under consideration, write $d\mu=\sigma\,dx$. Smoothness and homogeneity
of $F$ imply the local bounds
\[
    |\partial_xJ^*(x,\xi)|\leq C|\xi|, \qquad |d\log\sigma|\leq C,
\]
as well as the uniform ellipticity in \eqref{app:monotonicity}. These are the spatial coefficient assumptions in \cite[Section~4]{OhtaSturm2009}.

First, the equation for $u$ is a linear uniformly parabolic divergence-form equation with bounded right-hand side once its coefficients are evaluated along $du$. Interior parabolic estimates give a H\"older continuous
representative of $u$. A localized energy estimate with a test function of the form $\eta^2u_t$ gives $u_t\in L^2_{\mathrm{loc}}$ away from the initial time. This test is justified using the
time-discrete construction above, where convexity of $\mathcal E$ replaces differentiation of the energy. The additional term is bounded by
\[
    \left|\int_M\eta^2r(x,t,u)u_t\,d\mu\right| \leq\frac12\int_M\eta^2u_t^2\,d\mu +\frac12\int_M\eta^2r(x,t,u)^2\,d\mu.
\]
Apply next the local energy estimate to time difference quotients of $u$. Their diffusion matrices are averages of $\partial_\xi J^*$ along segments between the two gradients and are uniformly elliptic. Their reaction terms
consist of a bounded multiple of the difference quotient and a bounded term involving $r_t$. The argument in \cite[Proposition~4.5 and
Appendix~8.2]{OhtaSturm2009} therefore gives $u_t\in L^2_{\mathrm{loc}}(H^1_{\mathrm{loc}})$. Passing to the limit in the difference-quotient equation shows that $w=u_t$ satisfies
\begin{equation}\label{app:time-derivative}
    w_t=\Delta^{\nabla u}w +\bigl(a(\log u+1)+b\bigr)w +a_tu\log u+b_tu
\end{equation}
weakly. At $du=0$, one may choose any locally uniformly elliptic measurable extension of the reference matrix in this identity. Indeed, the Sobolev chain
rule gives $d_xu_t=\partial_t(d_xu)=0$ almost everywhere on that set, so the extension does not change the diffusion term. The zeroth-order coefficient and the inhomogeneous term in
\eqref{app:time-derivative} are locally bounded. Interior estimates for this linear equation give local boundedness and H\"older continuity of $u_t$.

Spatial difference quotients now give $u\in L^2_{\mathrm{loc}}(H^2_{\mathrm{loc}})$, as in \cite[Theorem~4.6 and Appendix~8.3]{OhtaSturm2009}. To see that the subsequent H\"older argument is also
unaffected, put $v=\partial_{x^i}u$. Differentiation in the weak sense gives
\[
    v_t=\operatorname{div}_\mu \bigl(\nabla^{\nabla u}v+H_i\bigr)+r_s(x,t,u)v+h_i,
\]
where
\[
\begin{split}
    H_i&=\partial_{x^i}J^*(x,du) +(\partial_{x^i}\log\sigma)J^*(x,du),\\
    h_i&=-(\partial_{x^i}\log\sigma)u_t +(\partial_{x^i}\log\sigma)r(x,t,u)+r_{x^i}(x,t,u).
\end{split}
\]
Here $\partial_{x^i}J^*$ differentiates the base variable with the covector fixed. Thus $|H_i|\leq C|du|$, while $h_i$ and $r_s$ are locally bounded. The
integrability iteration and H\"older estimate in \cite[Lemma~4.8 and Theorem~4.9]{OhtaSturm2009} apply to these equations. The bounded term $r_sv$ is absorbed in the same local energy
estimates by Gronwall's inequality. This proves the asserted H\"older continuity of $d_xu$. Finally, when $a$ and $b$ are smooth, the equation is smooth and uniformly parabolic on
every compact subset of $\{d_xu\neq0\}$, so interior regularity gives smoothness there. No smoothness assertion at critical points is needed.
\end{proof}

\end{document}